\documentclass[a4paper,11pt]{article}
\usepackage{amsmath,amsfonts,amssymb,amscd}
\usepackage{hyperref}
\newcommand\E{\mathrm{e}}

\newcommand\C{{\mathbb C}}
\renewcommand\P{{\mathbb P}}        
\newcommand\Q{{\mathbb Q}}  
\newcommand\R{{\mathbb R}}  
\newcommand\Z{{\mathbb Z}}   

\newcommand\germ[1]{{\mathfrak{#1}}}
\newcommand\h{{\germ h}}   
\newcommand\s{{\germ s}} 

\newcommand\phn{P_{\rm hn}}

\newcommand\proof{\noindent{\em Proof:\/ }}
\newcommand\qedsymbol{\hbox{\rlap{$\sqcap$}$\sqcup$}}
\newcommand\qed{\relax\ifmmode\else\unskip\quad\fi\qedsymbol}
 
\newtheorem{theorem}{Theorem}
\newtheorem{example}[theorem]{Example}
\newtheorem{definition}[theorem]{Definition}
\newtheorem{corollary}[theorem]{Corollary}
\newtheorem{lemma}[theorem]{Lemma}

\newtheorem{remark}[theorem]{Remark}
\renewcommand\qedsymbol{\hbox{\rlap{$\sqcap$}$\sqcup$}}
\renewcommand\qed{\relax\ifmmode\else\unskip\quad\fi\qedsymbol}
\newcommand\rarrow[3]{\smash{\mathop{\hbox to#3{\rightarrowfill}}\limits
^{\scriptstyle#1}_{\scriptstyle#2}}}
 
\newcommand\larrow[3]{\smash{\mathop{\hbox to#3{\leftarrowfill}}\limits
^{\scriptstyle#1}_{\scriptstyle#2}}}
 
\newcommand\darrow[3]{\llap{$\scriptstyle #1$}
\left\downarrow\vbox to#3{}\right.\rlap{$\scriptstyle #2$}}

\newcommand{\mapright}[1]{{\vbox{\ialign{
                                ##\crcr
    ${\scriptstyle\hfil\;\;#1\;\;\hfil}$\crcr
 \noalign{\kern2pt\nointerlineskip}
    \rightarrowfill\crcr}}\;}}

\newcommand{\mapleft}[1]{\mathop{\vbox{\ialign{
                                ##\crcr
    ${\scriptstyle\hfil\;\;#1\;\;\hfil}$\crcr
    \noalign{\kern2pt\nointerlineskip}
    \leftarrowfill\crcr}}\;}}

\newcommand{\mapdown}[1]{\Big\downarrow
         \rlap{$\vcenter{\hbox{$\scriptstyle#1$}}$}}

\newcommand\ve{{\scriptscriptstyle \vee}}

\newcommand\lset{\{}  
\newcommand\rset{\}}  
\newcommand\st{\ |\ }   
\newcommand\set[1]{\lset#1\rset} 
\newcommand\sett[2]{\lset #1 \st #2 \rset} 

\renewcommand\emptyset{\varnothing}
\renewcommand\setminus{-}
\newcommand\into{\hookrightarrow}

\newcommand\leftmapsto{\leftarrow\koinrel\hbox{\kern-3pt\raise1.5pt\hbox{\smash{$_\dashv$}}}}
\newcommand\longleftmapsto{\longleftarrow\koinrel\hbox{\kern-3pt\raise1.5pt\hbox{\smash{$_\dashv$}}}}
\newcommand\comp{\raise1pt\hbox{{$\scriptscriptstyle\circ$}}}

\newcommand\var{\mbox{\rm Var}}
\newcommand{\ii}{{\rm i}}
\newcommand\eh{\chi_{\mathrm{Hdg}}}
\newcommand\ehn{e_{\mathrm{Hdg}}}
\newcommand\eqdef{\mathop{:=}}
\newcommand\Gr{\mathop{\rm Gr}\nolimits}
\newcommand\Hom{\mathop{\rm Hom}\nolimits}
\newcommand\hs{\mathop{\h \s}\nolimits} 
\newcommand\fhs{\mathop{{\germ f} \h \s}\nolimits} 
\newcommand\Ker{\mathop{\rm Ker}\nolimits}

\title{Hodge number polynomials for nearby and vanishing cohomology}
\author{C.A.M.~Peters\\
Department of Mathematics,  University of Grenoble 
\\UMR 5582 CNRS-UJF, 38402-Saint-Martin d'H\`eres\\
France, email: chris.peters@ujf-grenoble.fr
\and 
J.H.M.~Steenbrink \\
Institute for Mathematics,Astrophysics and Particle Physics\\
Radboud University Nijmegen\\ 
Toernooiveld, NL-6525 ED Nijmegen\\
The Netherlands, email: J.Steenbrink@math.ru.nl}

\begin{document}
\maketitle
\begin{abstract}

\noindent We recall the construction of the Hodge character and we show,
using a result due to F.~Bittner,  that these can be constructed using
classical pure Hodge theory only, sideskipping Deligne's construction of
functorial mixed Hodge structures on the cohomology of algberaic varieties.
We then give a proof, based on the  weak factorization theorem, that the nearby
fibre is well-defined in the Grothendieck group of  varieties.  Finally,
various formulas for numerical invariants of Hodge-theoretic
origin are presented which are associated to  one-parameter degenerations of
projective manifolds. 
\end{abstract} 
{\small Keywords: Hodge structure, Hodge Euler polynomial, Hodge
number polynomial, nearby and vanishing cohomology}\\   {\small MSC2000
classification: 14D07, 32G20}

\section{Introduction}

The behaviour of the cohomology of a degenerating family of complex
projective manifolds has been intensively studied in the
nineteen-seventies by Clemens, Griffiths, Schmid and others. See
\cite{Gr} for a nice overview. Recently, the theory of motivic
integration, initiated by Kontsevich and developed by Denef and Loeser,
has given a new impetus to this topic. In particular, in the case of a
one-parameter degeneration it has produced an object $\psi_f$ in the
Grothendieck group of complex algebraic varieties, called the
\emph{motivic nearby fibre} \cite{B2}, which reflects the limit mixed
Hodge structure of the family in a certain sense. The purpose of this
paper is twofold. First, we prove that the motivic nearby fibre is
well-defined without using the theory of motivic integration. Instead
we use the Weak Factorization Theorem \cite{AKMW}. Second, we give a
survey of formulas containing numerical invariants of the limit mixed
Hodge structure, and in particular of the vanishing cohomology of an
isolated hypersurface singularity, without using the theory of mixed
Hodge structures or of variations of Hodge structure.

We hope that in this way  this interesting topic becomes accessible to
a wider audience. 
\section{Real Hodge structures}

A real Hodge structure on a finite dimensional real vector space $V$ consists of
a direct sum decomposition
\[V_\C=\bigoplus_{p,q\in \Z}V^{p,q}, \mbox{with $V^{p,q}=\overline{ V^{q,p}}$}  
\]
on its complexification $V_\C=V\otimes\C$. The corresponding 
\emph{Hodge filtration} is given by
\[
F^p(V)=\bigoplus_{r\ge p}V^{r,s}. 
\]
The numbers 
\[
h^{p,q}(V)\eqdef \dim V^{p,q}
\]
are the Hodge numbers  of the Hodge structure.  If for some integer $k$ we
have $h^{p,q}=0$ for all $(p,q)$ with $p+q\not=k$ the Hodge structure is pure
of weight $k$. Any real Hodge structure is the direct sum of pure Hodge
structures. The polynomial
\begin{eqnarray}
\phn(V)&=& \sum_{p,q\in\Z} 
h^{p,q}(V)u^pv^q\label{eqn:HEpoly} \\&=& \sum h^{p,k-p}(V)u^pv^{k-p}
\in\Z[u,v,u^{-1},v^{-1}]\nonumber 
\end{eqnarray}
is its associated  \emph{Hodge number polynomial}. \footnote{There are other
conventions in the litterature, for instance, some authors put a sign
$(-1)^{p+q}$ in front of the coefficient $h^{pq,}(V)$ of $u^pv^q$. } A classical
example  of a weight
$k$ Hodge real structure is furnished by the rank $k$ (singular) cohomology
group
$H^k(X)$ (with
$\R$-coefficients) of a compact K\"ahler manifold $X$.

Various multilinear algebra operations can be applied to Hodge structures as we
now explain. Suppose that $V$ and $W$ are two real vector spaces with a Hodge
structure of weight $k$ and $\ell$ respectively.  Then:
\begin{enumerate}
\vspace{-1ex}
\item $V\otimes W$ has a Hodge structure of weight $k+\ell$ given by
\[
F^p(V\otimes W)_\C =\sum_m
F^m(V_\C)\otimes F^{p-m}(W_\C)\subset V_\C\otimes_\C W_\C 
\]
and  with Hodge number polynomial given by
\begin{equation}
\phn(V\otimes W)= \phn(V)\phn(W). \label{eqn:HEforProds}
\end{equation}
\item 
On $\Hom(V,W)$ we have a Hodge structure of weight $\ell-k$:
\[F^p\Hom(V,W)_\C = \sett{f:V_\C\to
W_\C}{fF^n(V_\C)\subset F^{n+p}(W_\C)\quad\forall n}\]
with Hodge number polynomial
\begin{equation}
\phn(\Hom(V,W)(u,v)=\phn(V)(u^{-1},v^{-1})\phn(W)(u,v).
\label{eqn:HEforHoms}
\end{equation} 
In particular, taking $W=\R$ with $W_\C=W^{0,0}$ we get a Hodge structure of
weight $-k$ on the dual $V^\ve$ of $V$ with Hodge number polynomial 
\begin{equation}
\label{eqn:HEforDuals}\phn(V^\ve)(u,v)=\phn(V)(u^{-1},v^{-1}).
\end{equation} 
\end{enumerate}

The category $\h\s$ of real Hodge structures leads to a ring, the \emph{
Grothendieck ring} $K_0(\hs)$ which is  is the free group on the  isomorphism
classes
$[V]$ of real Hodge structures $V$  modulo the subgroup generated by
$[V]-[V']-[V'']$ 
where
\[
0\to V'\to V\to V''\to 0
\] 
is an exact sequence of pure Hodge structures and where the complexified maps
preserve the Hodge decompositions. Because the  Hodge number
polynomial(\ref{eqn:HEpoly}) is clearly additive  and    by
(\ref{eqn:HEforProds}) behaves well on products the  Hodge number polynomial
defines a ring homomorphism
\begin{equation*} 
\phn: \mathrm{K}_0(\hs) \to \Z[u,v,u^{-1},v^{-1}].
\end{equation*} 

As remarked before, pure Hodge structures of weight $k$ in algebraic geometry
arise as the (real) cohomology groups $H^k(X)$ of smooth complex projective
varieties. We combine these as follows: 
\begin{eqnarray} \label{eqn:chihodge}
\eh(X) &\eqdef & \sum (-1)^k [H^k(X)] \in \mathrm{K}_0(\hs); \\
\ehn(X) &\eqdef  & \sum (-1)^k\phn(H^k(X)) \in \Z[u,v,u^{-1},v^{-1}]
\label{eqn:hnpol} 
\end{eqnarray}
which we call the Hodge-Grothendieck character and the  Hodge-Euler
polynomial  of $X$ respectively.

\vspace{4mm}\noindent
Let us now recall the definition of the naive Grothendieck group
$\mathrm{K}_0(\var)$ of (complex) algebraic varieties. It is the quotient of the
free abelian group on isomorphism classes $[X]$ of algebraic varieties over $\C$
with the so-called \emph{scissor relations} $[X]=[X-Y]+[Y]$ for $Y\subset X$ a
closed subvariety. 
The cartesian product is compatible with the scissor relations and induces a
product structure on $\mathrm{K}_0(\var)$, making it into a ring.
There  is a nice set of generators and relations  for $\mathrm{K}_0(\var)$. To
explain this we first recall: 

\begin{lemma} \label{blow-up}
Suppose that $X$ is a smooth projective variety and $Y \subset X$ is a smooth
closed subvariety.  Let $\pi:Z \to X$ be the blowing-up with centre $Y$ and let
$E=\pi^{-1}(Y)$ be the exceptional  divisor. Then 
\begin{eqnarray*}
\eh(X)-\eh(Y)  &=& \eh(Z)-\eh(E);\\
\ehn(X)-\ehn(Y)  &=& \ehn(Z)-\ehn(E).
\end{eqnarray*}
\end{lemma}
\proof By \cite[p.~605]{G-H}
\begin{equation*}
0 \to H^k(X) \to H^k(Z) \oplus H^k(Y) \to H^k(E) \to 0
\end{equation*}
is exact. \qed  
 \begin{theorem}[\protect{\cite[Theorem 3.1]{B1}}] The group
$\mathrm{K}_0(\var)$ is isomorphic to the free abelian group generated by the
isomorphism classes of smooth complex projective varieties subject to the
relations $[\emptyset]=0$ and
$[Z]-[E]=[X]-[Y]$ where $X,Y,Z,E$ are as in Lemma \ref{blow-up}. 
\end{theorem}

It follows that   for every complex algebraic variety $X$ there
exist projective smooth varieties $X_1,\ldots,X_r,Y_1,\ldots,Y_s$ such
that 
\[
[X]=\sum_i[X_i]-\sum_j[Y_j] \mbox{ in } \mathrm{K}_0(\var) 
\]
and so, using Lemma~\ref{blow-up} we have:
\begin{corollary} The Hodge Euler character extends to a ring homomorphism
\[
\eh : K_0(\var) \to K_0(\h\s)
\]
and the Hodge number polynomial extends to  a ring homomorphism 
\[ 
\ehn: \mathrm{K}_0(\var) \to  \Z[u,v,u^{-1},v^{-1}] 
\]
\end{corollary}
\begin{remark} \rm
By  Deligne's theory
\cite{Del71}, \cite{Del74}  there is  a mixed Hodge structure  on the real
vector spaces $H^k(X)$. For our purposes, since we are working with real
coefficients, a mixed Hodge structure is just a real Hodge structure, i.e. a
direct sum of real Hodge structures of various weights, and so the Hodge
character and Hodge number polynomial are   defined for any real mixed Hodge
structure.  However, ordinary cohomology does not behave well with respect to
the scissor relation; we need compactly supported cohomology $H^k_c(X;\R)$. But
these also carry a Hodge structure and we have the following explicit
expression for the above characters.
\begin{eqnarray*} \eh(X) = \sum (-1)^k
[H^k_c(X))];\\ \ehn(X) = \sum (-1)^k \phn(H^k_c(X)).
\end{eqnarray*} 
\end{remark}
\begin{example} \rm 
\item{1)} Let $U$ be a smooth,
but not necessarily compact complex algebraic manifold.  Such a
manifold has   a \emph{good compactification} $X$, i.e. $X$ is a
compact complex algebraic manifold and $D=X\setminus U$ is a normal
crossing divisor, say $D=D_1\cup\cdots D_N$ with $D_j$  smooth and
irreducible. We introduce
\begin{eqnarray*} D_I&=& D_{i_1}\cap
D_{i_2}\cap\cdots \cap D_{i_m},\quad I=\set{i_1,\dots,i_m};\\a_I&:& D_I
\into X\end{eqnarray*}
and we set
\[
\begin{array}{rl}D(0)&=X;\\D(m)&=\coprod_{|I|=m}D_I,\quad
m=1,\dots,N;\\  a_m&=\coprod_{|I|=m}a_I:D(m)\to X.\end{array}
\]
Then each connected component of $D(m)$ is a  complex submanifold of $X$ of
codimension $m$.  Note that 
\[
[U] = \sum_m (-1)^m[D(m)] \in  \mathrm{K}_0(\var).\]
 Hence 
\begin{eqnarray*} \eh(U)& =& \sum_m (-1)^m
\eh(D(m));\\
\ehn(U)& = &\sum_m (-1)^m \ehn(D(m)).
\end{eqnarray*}
 \item{2)} If $X$ is compact  the construction of cubical 
hyperresolutions $(X_I)_{\emptyset\neq I \subset A}$ of $X$ from \cite{G-N-P-P}
leads to the  expression 
\[
[X]=\sum_{\emptyset\neq I \subset A}(-1)^{|I|-1} [X_I].
\]
and we find:
\begin{eqnarray*} \eh(X) & = &\sum_{\emptyset\neq I \subset
A}(-1)^{|I|-1} \eh(X_I); \\ 
\ehn(X) & = & \sum_{\emptyset\neq I
\subset A}(-1)^{|I|-1} \ehn(X_I). 
\end{eqnarray*}
\end{example}

The scissor-relations imply that the inclusion-exclusion principle can be
applied to a disjoint union $X$ of locally closed
subvarieties $X_1,\ldots,X_m$:  
\[  \eh(X)= \sum_{i=1}^m  \eh(X_i)
\]
and a similar expression holds for the Hodge Euler polynomials.

\begin{example}\rm  Let $T^n = (\C^\ast)^n$ be an  $n$-dimensional algebraic
torus. Then $\ehn(T^1) = uv-1$ so  $ \ehn(T^n)=(uv-1)^n$. Consider an
$n$-dimensional toric variety $X$. It is a  disjoint union of $T^n$-orbits.
Suppose that $X$ has $s_k$ orbits of dimension $k$. Then
\[
 \ehn(X) = \sum_{k=0}^n s_k \ehn(T^k) = \sum_{k=0}^n s_k(uv-1)^k.\]
If $X$ has a pure Hodge structure (e.g. if $X$ is compact  and has
only quotient singularities) then this formula determines the  Hodge
numbers of $X$.
\end{example}

\section{Nearby  and vanishing cohomology}

In this section we consider a relative stituation. We let $X$ be a complex
manifold, $\Delta\subset  \C$  the unit disk and 
$f:X \to \Delta$   a holomorphic map which is smooth over the punctured
disk $\Delta^*$. We say that $f$ is a \emph{one-parameter degeneration}.
Let us assume that  $E=\bigcup_{i\in I}E_i = f^{-1}(0)$ is a divisor
with strict normal crossings on $X$. We have the  \emph{specialization
diagram}
\[
\begin{matrix}X_\infty  & \mapright{k} & X & \mapleft{i} & E
\cr\mapdown{\tilde{f}} & & \mapdown{f} & & \mapdown{} \cr\h &
\mapright{\E} & \Delta  & \mapleft{} & \{0\} \cr\end{matrix}
\]
where $\h$ is the complex upper half plane, $\E(z):= \exp(2\pi\ii z)$ and
where   
\[
X_\infty\eqdef X\times_{\Delta^*} \germ{h}.
\] 
We let $e_i$ denote the multiplicity of $f$ along $E_i$ and choose a positive
integer multiple $e$ of all $e_i$. We let $\tilde{f}:\tilde{X}\to
 \Delta$ denote the normalization of the pull-back of $X$ under
the map $\mu_e:\Delta\to\Delta$ given by $\tau \mapsto \tau^e=t$. It fits into
a commutative diagram describing the \emph{$e$-th root} of $f$:
\[
\begin{matrix}
\tilde X & \mapright{\rho} & X\\
\darrow{\tilde f}{}{3ex}&& \darrow{f}{}{3ex}\\
\Delta&\mapright{\mu_e}& \Delta.
\end{matrix}
\]
We put 
\[
 D_i=\rho^{-1}(E_i)\, ,\,  D_J =
\bigcap_{i\in J} D_j,\,  D(m)=\coprod_{|J|=m}D_J. 
\]
Then $D_J\to E_J$ is a cyclic  cover of degree $\gcd(e_j\mid j\in J)$.
The maps $D_J \to E_J$ do not depend on the choice of the integer $e$ and so in
particular this is true for the varieties $D_J$ . See e.g.
\cite{Steen77a} or \cite{B2} for a detailed study of the geometry of
this situation. 
The special fibre $\tilde X_0=\tilde f^{-1}(0)$ is now a complex variety
equipped with the action of the cyclic group of order $e$. Let us introduce the
associated Grothendieck-group:

\begin{definition} \rm We let $\mathrm{K}^{\hat{\mu}}_0(\var)$ denote the
Grothendieck group of complex algebraic varieties with an action of a finite
order automorphism modulo the subgroup generated by expressions
$[\P(V)]-[\P^n\times X]$ where $V$ is a vector bundle of rank $n+1$ over $X$
with action which is linear over the action on $X$. 
See \cite[Sect.~2.2]{B2} for details.
\end{definition}
As a motivation, we should remark that in the ordinary Grothendieck group the
relation
$[\P(V)]=[\P^n\times X]$ follows from the fact that any algebraic vector bundle
is trivial over a Zariski-open subset and the fact that the scissor-relations
hold. The above relation extend to the case where one has a group action.

\begin{definition} \rm Suppose that the fibres of  $f$ are projective varieties.
Following \cite[Ch.~2]{B2} we define the \emph{motivic nearby fibre} of $f$ by
\begin{equation*} \psi_f := \sum_{m\geq 1} (-1)^{m-1}[D(m) \times
\P^{m-1}] \in \mathrm{K}^{\hat{\mu}}_0(\var)  
\end{equation*}
and the \emph{motivic vanishing fibre} by
\begin{equation}  \label{motivic}
\phi_f :=\psi_f-[E] \in \mathrm{K}^{\hat{\mu}}_0(\var)  
\end{equation}
\end{definition} 
\begin{remark} \rm If we let $E_I^0$ be the open subset of $E_I$ consisting of
points which are exactly on the $E_j$ with  $j\in I$,   $D_I^0$ the
corresponding subset of $D_I$ and $D^0(i)=\coprod_{|I|=i} D_I^0$. We have
$D_I=\coprod_{J\supset I} D^0_J$ and the scissor relations imply that 

\begin{eqnarray*}
 \sum_{i\ge 1} (-1)^{i-1} [ D(i)\times \P^{i-1}] &= &
\sum_{i\ge 1} (-1)^i \sum_{j\ge i} {j\choose i} [D^0(j)\times \P^{i-1} ]\\
&= &\sum_{j\ge 1} [D^0(j)] \times \sum_{i=1}^j (-1)^ {i-1} {j\choose i} [
\P^{i-1}]\\ & = &\sum_{j\ge 1} (-1)^{j-1}  [D^0(j)\times (\C^*)^{j-1}].
\end{eqnarray*}
This expression for the motivic nearby fibre has been used  in \cite{L}.
\end{remark}
 
The motivic nearby fibre turns out to be  a relative bimeromorphic invariant:

\begin{lemma} Suppose that $g:Y \to X$ is a
bimeromorphic proper map which is an isomorphism over $X-E$. Assume that
$(f\comp g)^{-1}(0)$ is a divisor with strict normal crossings. Then
\[
\psi_f = \psi_{fg}.\]
 \end{lemma} 
\proof In \cite{B2}, the proof relies on the theory of motivic integration
\cite{D-L}. We give a different proof, based on the \emph{weak factorization
theorem} \cite{AKMW}. This theorem reduces the problem to the following
situation: $g$ is the blowing-up of $X$ in a connected submanifold $Z\subset E$,
with the following property. Let $A\subset I$ be those indices $i$ for
which $Z\subset E_i$. Then $Z$ intersects the divisor $\bigcup_{i\not\in A}E_i$
transversely, hence $Z\cap \bigcup_{i\not\in A}E_i$ is a divisor with normal
crossings in $Z$.

We fix the following notation. We let $L:=[\mathbb{A}^1]$ and
$P_m:=[\P^m]=\sum_{i=0}^mL^i$. For $J \subset I$ we let $j=|J|$.
So, using the product structure in $K_0(\var)$, we have 
\[ \psi_f = \sum_{\emptyset \neq J \subset I}
(-1)^{j-1}[D_J]P_{j-1}. 
\]
Let $E'= \bigcup_{i\in I'}E'_i$ be the zero fibre of $fg$. We have $I=I'\cup
\{\ast\}$ where $E'_*$ is the  exceptional divisor of $g$ and $E'_i$ is the
proper transform in $Y$ of $E_i$.  Form $\rho':D'\to E'$, the associated
ramified cyclic covering. For $J\subset I$ we let $J'=J\cup\{\ast\}$. Note that
$E'_\ast$ has multiplicity equal to 
$\sum_{i\in
A}e_i$. 
Without loss of generality we may assume that $e$ is also a
multiple of this integer. We have two kinds of $j'$-uple intersections
$D_{J'}$: those which only contain $D'_j$, $j\not=\ast$ and those which contain
$D'_*$. So,
\[ \psi_{fg} = [D'_\ast]
+\sum_{\emptyset \neq J \subset
I}(-1)^{j-1}\left([D'_J]P_{j-1}-[D'_{J'}]P_j\right). 
\]
%
%
We are going to
calculate the difference between $\psi_{fg}$ and $\psi_f$. Let $B\subset
I-A$. Let $c=\mathrm{codim}(Z,X)$. Then for all $K\subset A$ we have
that $\mathrm{codim}(Z\cap E_{K\cup B}, E_{K\cup B})= c-k$, and
$D'_{K\cup B}$ is the blowing up of $D_{K\cup B}$ with centre $Z\cap
E_{K\cup B} = Z\cap E_B=:Z_B$ and exceptional divisor $D'_{K\cup B'}$.
Hence if we let $W_B={\rho'}^{-1}(Z_B)$ we have 
\[
[D'_{K\cup B}] = [D_{K\cup B}] +  [D'_{K\cup B'}]-[W_B] = [D_{K\cup B}] +
[W_B](P_{c-k-1}-1).
\]
Hence 
\[
\psi_{fg}-\psi_f = \sum_B c_B[W_B]
\]
for suitable coefficients $c_B$.

Note that $Z_\emptyset=Z$ and that $D'_*=W_\emptyset \times P_{c-1}$ and so if
$B=\emptyset$ we get 
\begin{eqnarray*} c_\emptyset & = & P_{c-1}+
\sum_{\emptyset \neq K \subset A} (-1)^{k-1}\left( (P_{c-k-1}
-1)P_{k-1}-P_{c-k-1}P_k \right)\\ & = & P_{c-1}+ \sum_{\emptyset \neq K
\subset A} (-1)^k (L^kP_{c-k-1}+P_{k-1}) = P_{c-1}\sum_{K\subset A}
(-1)^k = 0. 
\end{eqnarray*} 
In case  $B\neq \emptyset$ we get
\begin{eqnarray*} c_B & = & \sum_{K \subset A} (-1)^{k+b-1}\left(
(P_{c-k-1} -1)P_{k+b-1}-P_{c-k-1}P_{k+b} \right)\\ &= &\sum_{K\subset
A} (-1)^{k+b} (L^{k+b}P_{c-k-1}+P_{k+b-1}) = P_{c+b-1}\sum_{K\subset A}
(-1)^{b+k} = 0.  \quad \qed
\end{eqnarray*}

Let us now pass to the nearby fibre in the Hodge theoretic sense.
In \cite{Sch73} and \cite{Steen76} a mixed Hodge
structure on $H^k(X_\infty)$ was constructed; its weight filtration is
the \emph{monodromy weight filtration} which we now explain. 
The loop winding once counterclockwise around the origin gives a generator
of
$\pi(\Delta^*,*)$ ,$*\in \Delta^*$. Its action on the fibre $X_*$
over $*$ is well defined up to homotopy and on $H^k(X_*)\simeq H^k(X_\infty)$
it defines the \emph{monodromy automorphism} $T$. Let $N = \log T_u$
be the logarithm of the unipotent part in the Jordan decomposition of $T$.
Then $W$ is the unique increasing filtration on
$H^k(X_\infty)$ such that $N(W_j)\subset W_{j-2}$ and $N^j:\Gr^W_{k+j}
\to \Gr^W_{k-j}$ is an isomorphism for all $j\geq 0$. In fact, from
%
%
\cite[Lemma 6.4]{Sch73} we deduce:
\begin{lemma}  \label{jordan} There is a Lefschetz-type decomposition 
\[
\Gr^W H^k(X_\infty)= \bigoplus_{\ell=0}^k \bigoplus_{r=0}^{\ell} N^r P_{k+\ell}, 
\]
where $P_{k+\ell}$ is pure of weight $k+\ell$. The endomorphism $N$ has $\dim
P_{\ell+m-1}$ Jordan blocs of size $m$.
\end{lemma}
The Hodge filtration is constructed in \cite{Sch73} as a limit of
the Hodge filtrations on nearby smooth fibres in a certain sense, and in
\cite{Steen76} using the relative logarithmic de Rham complex. With $X_t$   a
smooth fibre of $f$ it follows that
\[
\dim F^m H^k(X_t)= \dim F^m H^k(X_\infty)
\]
and hence
\begin{equation} \label{eqn:hodgefilt}
\ehn(X_\infty)_{|v=1} = \ehn(X_t)_{|v=1}.
\end{equation}
We next remark that the nearby cycle sheaf $i^*Rk_*k^*\R_X$ can be used  to put
a    Hodge structure on the cohomology groups $H^k(X_\infty)$, while the
hypercohomology of the vanishing cycle sheaf $\phi_f=\mbox{\rm Cone}[i^*\R_X\to
\psi_f\R_X]$ (it is a \emph{complex} of sheaves of real vector spaces on
$E$) likewise admits a Hodge structure. In fact, the spectral
sequence of \cite[Cor.~4.20]{Steen76} shows that
\[
\eh(X_\infty) =  \eh(\psi_f)
\]
and then the definition shows that
\[
\eh(X_\infty) - \eh(E)  =  \eh(\phi_f) .
\]
In fact this formula motivates the nomenclature ``motivic nearby fibre'' and
``motivic vanishing cycle''. 

As a concluding remark, the semisimple part $T_s$ of the monodromy
is an automorphism of the mixed Hodge structure on $H^k(X_\infty)$. 

We can use these remarks to deduce information about the Hodge numbers on
$H^k(X_\infty)$ from information about the geometry of the central fibre as we
shall illustrate now.

\begin{example}\rm 

\item{1)} Let $F,L_1,\ldots,L_d \in\C[X_0,X_1,X_2]$ be 
homogeneous forms with $\deg F=d$ and $\deg L_i=1$  for $i=1,\ldots,d$,
such that $F\cdot L_1\cdots L_d=0$ defines a reduced divisor with normal
 crossings on $\P^2(\C)$. We consider the space \[X=\{([x_0,x_1,x_2],t)
\in \P^2\times
\Delta\mid\prod_{i=1}^dL_i(x_0,x_1,x_2)+tF(x_0,x_1,x_2)=0\}\] where
$\Delta$ is a small disk around $0\in\C$. Then $X$ is smooth and the
map$f:X\to \Delta$ given by the  projection to the second factor has as
its zero fibre the union $E_1\cup\cdots\cup E_d$ of the lines  $E_i$
with equation$L_i=0$. These lines are in general position and have
multiplicity one. We obtain 
\[ \psi_f = (d-{d \choose 2})[\P^1] \]
 so
\[
\ehn(\psi_f )=(1-{d-1 \choose 2})(1+uv)
\]
and substituting $v=1$ in this formula we get the formula $g={d-1 \choose 2}$
for the genus of a smooth  plane curve of degree $d$. The monodromy on
$H^1(X_\infty )$ has $g$ Jordan blocks of size $2$, so is ``maximally 
unipotent''.
  
\item{2)}
If we consider a similar example,  but
replace $\P^2$ by $\P^3$ and curves by surfaces, lines by planes, then
the  space $X$ will not be smooth but has ordinary double points at the
points of the zero fibre where two of the  planes meet the surface
$F=0$. There are $d{d\choose 2}$ of such points, $d$ on each line of
intersection. If we blow these up, we obtain a family $f:\tilde{X}
\to\Delta$ whose zero fibre $D=E\cup F$ is the union  of components
$E_i,\ i=1,\ldots,d$ which are copies of $\P^2$ blown up in $d(d-1)$
points, and components $F_j, j=1,\ldots,d{d\choose 2}$ which are copies
of $\P^1\times\P^1$.  Thus 
\[
\ehn(D(1))= d(1+(d^2-d+1)uv+u^2v^2)+ d{d\choose2}(1+uv)^2 .
\] 
The double point locus $D(2)$ consists of  the
${d \choose 2}$ lines of intersections of the $E_i$ together with the
$d^2(d-1)$ exceptional lines in the $E_i$. So 
\[
\ehn(D(2))=d(d-1)(d+\frac 12)(1+uv) .
\]
Finally $D(3)$ consists of the ${d \choose 3}$ intersection points of the $E_i$
together with one point on each  component $F_j$, so 
\[
\ehn (D(3))={d\choose 3}+ d{d\choose 2}= \frac 13 d(d-1)(2d-1) \]
We get 
\[
\ehn(\psi_f)=\left({d-1\choose 3}+1\right)(1+u^2v^2)+ \frac 13 d(2d^2-6d+7)uv
\]
in accordance
with the Hodge numbers for a smooth degree $d$ surface:
\[
h^{2,0}=h^{0,2}={d-1\choose 3}, \quad h^{1,1}=\frac 13 d(2d^2-6d+7)
\]
The monodromy on $H^2(X_\infty )$ has $\displaystyle {d-1\choose 3}$
Jordan blocs of size $3$ and 
%
%
$\frac 12 d^3-  d^2+\frac 12 d +1$ 
blocks of size
$1$. 

\item{3)} Consider a similar  smoothing of
the union of two transverse quadrics in $\P^3$. The generic fibre is a
smooth  K3-surface and after blowing up the 16 double points of the
total space we obtain the following special fibre:
\begin{enumerate}
\item $E(1)$ has  two components which are blowings up
of $\P^1\times\P^1$ in 16 points, and 16 components isomorphic to
$\P^1\times\P^1$; hence $\ehn(E(1))=18(1+uv)^2 + 32uv$.
\item $E(2)$ consists of the 32  exceptional lines together with the strict
transform of the  intersection of the two quadrics, which is an elliptic
curve; hence   $\ehn(E(2))=33(1+uv)-u-v$;\item $E(3)$ consists of 16
points:  one point on each exceptional $\P^1\times\P^1$, so
$\ehn(E(3))=16$.
\end{enumerate}
We get 
\[
\ehn(X_\infty)=1+u+v+18uv+u^2v+uv^2+u^2v^2
\]
Putting $v=1$ we get $2+20u+2u^2$, in agreement with the Hodge numbers
$(1,20,1)$ on the $H^2$ of a K3-surface.  The monodromy has two Jordan blocs of
size $2$ and $18$ blocs of size $1$.
\end{example}

\section{Equivariant Hodge number polynomials} 
We have seen that the mixed Hodge structure on the cohomology of the nearby
fibre of a one-parameter degeneration comes with an automorphism of
finite order. This leads us to consider the category
$\hs_\R^{\hat{\mu}}$ of pairs $(H,\gamma)$ consisting of a  real Hodge
structure (i.e. direct sum of pure real Hodge structures of possibly
different weights) $H$ and an automorphism $\gamma$ of finite order of
this  Hodge structure.

We are going to consider a kind of tensor product of two such objects,
which we call \emph{convolution} (see \cite{SchS85}, where this
operation was defined for mixed Hodge structures and called
\emph{join}). We will explain this by settling an equivalence of
categories between $\hs_\R^{\hat{\mu}}$ and a category $\fhs$ of
so-called \emph{fractional  Hodge structures}. (called \emph{Hodge
structures with fractional weights} in \cite{L}; it is however not the
weights which are fractional, but the indices of the Hodge filtration!

\begin{definition} \rm (See \cite{L}). A \emph{fractional Hodge
structure} of weight $k$ is a real vector space $H$ of finite
dimension, equipped with a decomposition 
\[H_\C = \bigoplus_{a+b=k}
H^{a,b}\]
 where $a,b\in \Q$, such that $H^{b,a} = \overline{H^{a,b}}$.
A \emph{fractional Hodge structure} is defined as a direct sum of pure
fractional Hodge structures of possibly different weights.
\end{definition}

\begin{lemma} We have an equivalence of categories  $G: \hs^{\hat{\mu}}
\to \fhs$. \end{lemma} \proof Let $(H,\gamma)$ be an object of
$\hs^{\hat{\mu}}$ pure of weight $k$. We define $H_a =
\Ker(\gamma-\exp(2\pi\ii a);H_\C)$ for $0\leq a <1$ and for $0<a<1$ put
%
%
\[
\tilde{H}^{p+a,k+1-a}= H_a^{p,k-p}, \quad \tilde{H}^{p,k-p}=H_0^{p,k-p}\]
 This transforms $(H,\gamma)$ into a
direct sum $\tilde{H} =:G(H,\gamma)$ of fractional Hodge structures of
weights $k$ and $k+1$ respectively. Conversely, for a fractional Hodge
structure $\tilde{H}$ of weight $k$ one has a unique automorphism
$\gamma$ of finite order which is multiplication by $\exp(2\pi\ii b)$
on $\tilde{H}^{b,k-b}$.

\vspace{4mm}\noindent Note that this equivalence of categories does not
preserve tensor products! Hence it makes sense to make
\begin{definition}
\rm The \emph{convolution} $(H',\gamma')\ast
(H",\gamma")$ of two objects in $\hs^{\hat{\mu}}$ is the object
corresponding to the tensor product of their images in $\fhs$:
\[ 
G\left( (H',\gamma') \ast (H",\gamma") \right) = G(H',\gamma') 
\otimes G(H",\gamma").
\]
\end{definition} 
Note that the Hodge number
polynomial map $\phn$ extends to a ring homomorphism 
\[
\phn^{\hat{\mu}}: K_0(\fhs) \to R:=\lim_{\leftarrow}\Z[u^{\frac
1n},v^{\frac 1n},u^{-1},v^{-1}]
\]
 We denote its composition with the
functor $G$ by the same symbol. Hence
\[ \phn^{\hat{\mu}}:
K_0^{\hat{\mu}}(\hs) \to R \]
 transforms convolutions into products. We
equally have an equivariant Hodge-Grothendieck character
\[
\chi_\mathrm{Hdg}^{\hat{\mu}}: K_0^{\hat{\mu}}(\var) \to
K_0^{\hat{\mu}}(\h\s)\]
 and an equivariant Hodge-Euler characteristic
\[
\ehn^{\hat{\mu}}:  K_0^{\hat{\mu}}(\var) \to  R.
\]
 
Let $f:X \to \C$ be a projective morphism where $X$ is smooth, of
relative dimension $n$, with a single isolated critical point $x$ such
that $f(x)=0$. Construct $\psi_f$ by replacing the zero fibre $X_0$ by
a divisor with normal crossings as above.   Then the Milnor fibre of $F$ of $f$
at $x$ has the homotopy type of a wedge of spheres of dimension $n$.
Its cohomology is also equipped with a mixed Hodge structure. Recalling
 definition~\ref{motivic} of the motivic vanishing fibre, it can be shown that
\[
\phn^{\hat{\mu}}(\tilde{H}^n(F)) = (-1)^n\ehn^{\hat{\mu}}(\phi_f).
\]
Write $\phn^{\hat{\mu}}(\tilde{H}^n(F)) = \sum_{\alpha\in\Q,\ w\in\Z}
m(\alpha,w)u^\alpha v^{w-\alpha}$. In the literature several numerical
invariants have been attached to the singularity $f:(X,x)\to (\C,0)$. These
are all related to the numbers $m(\alpha,w)$ as follows: 
\begin{enumerate}
\item The \emph{characteristic pairs} \cite[Sect.~5]{Steen77a}.
\[
\mathrm{Chp}(f,x) = \sum_{\alpha,w} m(\alpha,w)\cdot (n-\alpha,w).
\]
\item The \emph{spectral pairs} \cite{N-S}: 
\[\mathrm{Spp}(f,x) =
\sum_{\alpha\not\in\Z, w} m(\alpha,w)\cdot (\alpha,w) +
\sum_{\alpha\in\Z,w}m(\alpha,w)\cdot (\alpha,w+1).\]
 
\item The singularity spectrum in Saito's sense \cite{Sa83}:
\[
\mathrm{Sp}_\mathrm{Sa}(f,x) =
\phn^{\hat{\mu}}(\tilde{H}^n(F))(t,1).\]

\item The singularity spectrum
in Varchenko's sense \cite{Var81}: 
\[
\mathrm{Sp}_\mathrm{V}(f,x) =
t^{-1}\phn^{\hat{\mu}}(\tilde{H}^n(F))(t,1).
\]
 \end{enumerate}

Note that the object $\phi_f$ depends only on the germ of
$f$ at the critical point $x$. Let us as a final remark  rephrase  the
original   Thom-Sebastiani theorem (i.e. for the case of isolated
singularities):
\begin{theorem} Consider holomorphic germs $f:(\C^{n+1},0)\to (\C,0)$ and
$g:(\C^{m+1},0) \to (\C,0)$ with isolated singularity. Then the germ
$f\oplus g: (\C^{n+1}\times\C^{m+1},(0,0)) \to \C$ with $(f\oplus
g)(x,y):=f(x)+g(y)$ has also an isolated singularity, and
\[
\eh^{\hat{\mu}}( \phi_{f\oplus g}) = -\eh^{\hat{\mu}}(\phi_f) \ast
\eh^{\hat{\mu}}(\phi_g )
\]
 so 
\[
\ehn^{\hat{\mu}}( \phi_{f\oplus g}) =
-\ehn^{\hat{\mu}}(\phi_f) \cdot \ehn^{\hat{\mu}}(\phi_g )\in R.
\]
\end{theorem} 
\begin{remark}\rm This theorem has been largely
generalized, for functions with arbitrary singularities, and even on
the level of motives. Denef and Loeser \cite{D-L} defined a convolution
product for Chow motives, and Looijenga \cite{L} defined one on
$\mathcal{M}^{\hat{\mu}}=K_0^{\hat{\mu}}(\var)[L^{-1}]$, both with the
property that the Hodge Euler polynomial commutes with convolution and
that the Thom-Sebastiani property holds already on the level of
varieties/motives. 
\end{remark}

\end{document}